\documentclass[12pt]{article}
\usepackage{amsmath, amssymb, graphics, graphicx, color}
\usepackage{dsfont}
\newif\ifextended
\ifx\Extended\Kaba\else\extendedtrue\fi
\newif\ifJapanese
\newif\iftesting
\testingfalse

\usepackage{theorem}
\newcommand{\bbd}[1]{{\mathbb{#1}}}
\theorembodyfont{\ifJapanese\rm\else\it\fi}
\newcount\minute	
\newcount\hour		
\newcount\hourMins  
\ifJapanese
\def\today%
{
  \the\year 年\,\zeroPadTwo{\the\month}月\,\zeroPadTwo{\the\day}日%
}
\fi
\def\now%
{
  \minute=\time    
  \hour=\time \divide \hour by 60 
  \hourMins=\hour \multiply\hourMins by 60
  \advance\minute by -\hourMins 
  \zeroPadTwo{\the\hour}:\zeroPadTwo{\the\minute}%
}
\def\zeroPadTwo#1%
{
  \ifnum #1<10 0\fi    
  #1
}

\renewcommand{\baselinestretch}{1.12}

\let\Label\label%
\iftesting
\def\label#1{\marginpar{{\renewcommand{\baselinestretch}{0.6}\tiny 
			#1}}\Label{#1}\ignorespaces}%
\fi
\newcounter{frml}[section]
\newcounter{frmla}[section]
\def\thefrml{{\arabic{section}.\arabic{frml}}}
\def\thefrmla{{\arabic{section}.\arabic{frmla}a}}
\def\frmlabel#1{\refstepcounter{frml}{\def\baka{#1}\ifx\baka\empty\else\Label{#1}\fi}%
{\rm({\thefrml})\hfill\hfill\hfill}}
\def\frmlabela#1{\refstepcounter{frmla}{\def\baka{#1}\ifx\baka\empty\else\Label{#1}\fi}%
{\rm({\thefrmla})\hfill\hfill\hfill}}
\def\xitem[#1]{\item[\frmlabel{#1}]\mbox{}%
	\iftesting\marginpar{{\renewcommand{%
				\baselinestretch}{0.6}\tiny#1}}\fi\ignorespaces}
\def\xitemq[#1]{\item[\frmlabel{#1}]\mbox{}%
	\ignorespaces}
\def\xitema[#1]{\item[\frmlabela{#1}]\mbox{}%
	\iftesting\marginpar{{\renewcommand{%
				\baselinestretch}{0.6}\tiny#1}}\fi\ignorespaces}
\def\xitemsub[#1]#2{\item[\frmlabel{#1}$_{#2}$]\mbox{}%
	\iftesting\marginpar{{\renewcommand{%
				\baselinestretch}{0.6}\tiny#1}}\fi\ignorespaces}
\def\xxitem[#1][#2]{\item[(\ref{#1}{\makebox[1.4ex][c]{#2}})]\mbox{}%
	\iftesting\marginpar{{\renewcommand{%
				\baselinestretch}{0.6}\tiny\{#1\}\{#2\}}}\fi\ignorespaces}
\def\xitemof#1{{\rm({\ref{#1}})}}

\newenvironment{xitemize}{\begin{list}{}{
                        \parsep=0.5\smallskipamount%
			\itemindent=-0.4ex%
			\itemsep=0.5\smallskipamount\leftmargin=4ex
			\labelwidth=3ex
                        \labelsep=0.7ex                        
}}%
							 {\end{list}}
\def\assert#1{\noindent\makebox[4.8ex][r]{\rm(\makebox[1.2ex][c]{#1})\ \ }%
	\ignorespaces}
\def\assertof#1{{\rm (#1)}}

\newenvironment{proof}{{\it Proof.\ \ }}{\mbox{}\hfill$\Box$}
\newtheorem{Def}{\ifJapanese{\bf 定義}\else {\bf Definition}\fi}[section]
\newtheorem{definition}[Def]{\ifJapanese{\bf 定義}\else {\bf Definition}\fi}
\newtheorem{Thm}{\ifJapanese{\bf 定理}\else {\bf Theorem}\fi}[section]
\newtheorem{theorem}[Thm]{\ifJapanese{\bf 定理}\else {\bf Theorem}\fi}

\newtheorem{example}[Thm]{\ifJapanese{\bf 例}\else {\bf Example}\fi}

\newtheorem{problem}[Thm]{\ifJapanese{\bf 未解決問題}\else{\bf Problem}\fi}

\newtheorem{lemma}[Thm]{\ifJapanese{\bf 補題}\else{\bf Lemma}\fi}

\newtheorem{corollary}[Thm]{\ifJapanese{\bf 系}\else{\bf Corollary}\fi}

\newtheorem{Claim}{{\bf Claim}}[Thm]
\newtheorem{claim}[Claim]{{\bf Claim}}

\newcommand{\prf}{\ifJapanese{\bf 証明．\ }\ignorespaces\else{\bf 
		Proof.\ \ }\ignorespaces\fi}

\newcommand{\prfofClaim}{\raisebox{-.4ex}{\Large $\vdash$\ \ }}

\newcommand{\Thmof}[1]{\ifJapanese{定理\,\ref{#1}}\else{Theorem~\ref{#1}}\fi}

\newcommand{\bfThmof}[1]{\ifJapanese{\bf 定理\,\ref{#1}}\else{\bf Theorem~\ref{#1}}\fi}
\newcommand{\itThmof}[1]{\ifJapanese{\it 定理\,\ref{#1}}\else{\it Theorem~\ref{#1}}\fi}
\newcommand{\Lemmaof}[1]{\ifJapanese{補題\,\ref{#1}}\else{Lemma~\ref{#1}}\fi}

\newcommand{\Corof}[1]{\ifJapanese{系\,\ref{#1}}\else{Corollary~\ref{#1}}\fi}

\newcommand{\Thmabove}{{\ifJapanese 定理\else Theorem\fi\ \number\theThm}}

\newsavebox{\qedbox}\sbox{\qedbox}{
{\unitlength=0.05mm \begin{picture}(40,60)
\put(0,0){\framebox(30,44)[cc]{}}
\put(30,-7){\rule{7\unitlength}{44\unitlength}}
\put(10,-7){\rule{27\unitlength}{7\unitlength}}
\end{picture}}}
\newcommand{\qed}{\mbox{}\hfill\usebox{\qedbox}}
\newcommand{\smallqed}%
{\mbox{}\smallskip\hfill\raisebox{-.4ex}{\Large $\dashv$}}
\newcommand{\qedof}[1]%
{\mbox{} \hspace*{\fill}{\usebox{\qedbox}{\tiny~(#1)}}}
\newcommand{\Qedof}[1]%
{\mbox{} \hspace*{\fill}{\usebox{\qedbox}%
{\tiny~(#1~\number\theThm)}}}

\newcommand{\qedskip}{\medskip}

\newcommand{\qedofClaim}%
{\mbox{}\hfill\raisebox{-.4ex}{\Large $\dashv$ }\nolinebreak%
\mbox{\tiny~(Claim~\number\theClaim)}}
\newcommand{\qedofSubclaim}%
{\mbox{}\hfill\raisebox{-.4ex}{\Large $\dashv$ }\nolinebreak%
\mbox{\tiny~(Subclaim~\number\theSubclaim)}}
\newcommand{\cardof}[1]{\mathopen{|\,}#1\mathclose{\,|}}

\newcommand{\setof}[2]{\{#1\,:\,#2\}}
\newcommand{\ssetof}[1]{\{#1\}}

\newcommand{\subseteqand}[1]{\mathrel{\mathop{\subseteq}%
		\limits_{\scriptscriptstyle\hbox to 14pt{$\scriptscriptstyle #1$\hss}}}}

\newcommand{\mapping}[3]{#1:#2\rightarrow #3}

\newcommand{\fnsp}[2]{\mbox{}^{{#1}\hspace{-0.02em}}#2}
\newcommand{\imageof}{{}^{\,{\prime}{\prime}}}
\newcommand{\seqof}[2]{\langle#1\,:\,#2\rangle}
\newcommand{\pairof}[1]{\langle#1\rangle}

\newcommand{\forces}[2]{\,\|\hspace{-.35ex}\mbox{\sf--}_{\,#1\,}%
\mbox{\rm``}\,#2\,\mbox{\rm''}}

\newcommand{\restr}{\restriction}

\newcommand{\Fn}{{\rm Fn}}
\newcommand{\Cohen}[1]{{\calC}_{#1}}
\newcommand{\random}[1]{{\calR}_{#1}}
\newcommand{\flr}[1]{\mathopen{\lceil\hspace{-0.1ex}}#1\mathclose{\hspace{0.1ex}\rceil}}

\newcommand{\dom}{\mathop{\rm dom}}

\newcommand{\concat}{\mathop{{}^{\frown}}}

\newcommand{\calND}{{\mathcal N\hspace{-0.2ex}\mathcal D}}
\newcommand{\calNDN}{{\mathcal N\hspace{-0.2ex}\mathcal D\hspace{-0.2ex}\mathcal N}}
\newcommand{\poP}{\bbd{P}}
\newcommand{\poQ}{\bbd{Q}}

\newcommand{\cov}{{\sf cov}}

\newcommand{\non}{{\sf non}}

\newcommand{\gma}{\mathfrak{a}}
\newcommand{\gmb}{\mathfrak{b}}
\newcommand{\gmd}{\mathfrak{d}}
\newcommand{\gms}{\mathfrak{s}}
\newcommand{\gmo}{\mathfrak{o}}

\newcommand{\continuum}{\mathfrak{c}}

\newcommand{\calA}{{\mathcal A}}
\newcommand{\calB}{{\mathcal B}}
\newcommand{\calC}{{\mathcal C}}
\newcommand{\calD}{{\mathcal D}}

\newcommand{\calF}{{\mathcal F}}
\newcommand{\calG}{{\mathcal G}}
\newcommand{\calH}{{\mathcal H}}
\newcommand{\calI}{{\mathcal I}}

\newcommand{\calM}{{\mathcal M}}
\newcommand{\calN}{{\mathcal N}}
\newcommand{\calO}{{\mathcal O}}
\newcommand{\calP}{{\mathcal P}}
\newcommand{\calR}{{\mathcal R}}
\newcommand{\calS}{{\mathcal S}}

\newcommand{\calX}{{\mathcal X}}

\newcommand{\ZFC}{{\sf ZFC}}
\newcommand{\CH}{{\sf CH}}

\newcommand{\MA}{{\sf MA}}

\newcommand{\st}{such that}
\newcommand{\wrt}{with respect to}

\newcommand{\wolog}{without loss of generality}

\newcommand{\po}{poset}
\newcommand{\pos}{posets}

\newcommand{\utildeT}[1]{%
	\hbox to 0pt{$\mathop{#1}\limits_{\raise0.25ex\hbox{$\scriptstyle\sim$}}$\hss}%
		\relax\phantom{\underline{#1}}}
\newcommand{\utildeS}[1]{%
	\hbox to 0pt{\smash{$\mathop{\scriptstyle #1}\limits_{%
				\raisebox{0.6ex}[0pt]{$\scriptscriptstyle\sim$}}$}\hss}%
		\relax\phantom{\mathord{{#1}_{\rule[-0.6ex]{0pt}{1pt}}}}}
\newcommand{\utildeSS}[1]{%
	\hbox to 0pt{$\mathop{\scriptscriptstyle #1}%
		\limits_{\scriptscriptstyle\sim}$\hss}%
		\relax\phantom{\underline{#1}}}
\newcommand{\utilde}[1]{%
	\mathchoice{\utildeT{#1}}{\utildeT{#1}}{\utildeS{#1}}{\utildeSS{#1}}}

{\end{minipage}\end{trivlist}}
\makeatletter
\def\blfootnote{\xdef\@thefnmark{}\@footnotetext}
\makeatother
\newcommand{\xthanks}[1]{{\def\at{\quad}\def\and{\\{}\qquad}\blfootnote{#1}}}
\newcommand{\institute}[1]{{\def\at{\quad}\def\and{\\{}\qquad}\blfootnote{#1}}}
\newcommand{\keywords}[1]{\blfootnote{Keywords: #1}}
\begin{document}
\title{How to drive our families mad}
%
\author{Saka\'e Fuchino \and Stefan Geschke \and Osvaldo Guzman \and Lajos Soukup
}
\maketitle
\date

\xthanks{The first author was partially supported by Chubu 
	University grant 16IS55A,  as well as Grant-in-Aid for 
	Scientific Research (C) 19540152 and Grant-in-Aid for Exploratory 
  Research No.\ 26610040 of Japan Society for the Promotion of 
	Science. 
  The second author was supported by Centers of 
	Excellence grant from the European Union.  
  The third author was supported by CONACyT scholarship 420090. 
  The fourth author was supported by Bolyai Grant. The research of this paper 
	began when the first and second authors visited Alfr\'ed R\'enyi  
	Institute in Budapest in 2002. The research was then resumed when the first 
	author visited Centre de Recerca Matem\`atica in Barcelona and 
	the fourth author Barcelona University at the same time in 2006.   
	The first and fourth authors would like to thank 
	Joan Bagaria and Juan-Carlos Mart\'\i nez for the arrangement of the 
	visit and their hospitality during the stay in Barcelona. The authors 
	also would like to thank Kenneth Kunen for allowing them to include his unpublished 
	results and Andreas Blass as well as the referee of the paper for reading carefully the 
  manuscript and giving many valuable suggestions. }






\begin{abstract}
	Given a family $\calF$ of pairwise almost disjoint (ad) sets on a 
	countable set $S$,  
	we study maximal almost disjoint 
	(mad) families $\tilde{\calF}$ extending $\calF$.

	We define $\gma^+(\calF)$ to be the minimal possible cardinality 
	of $\tilde{\calF}\setminus \cal F$ for such $\tilde{\calF}$ and
	$\gma^+(\kappa)=\max\setof{\gma^+(\calF)}{\cardof{\calF}\leq\kappa}$. 
	We show that all infinite cardinals less than or equal to the continuum
	$\continuum$ can be 
	represented as $\gma^+(\calF)$ for some ad $\calF$ (\Thmof{aplus-all}) 
	and that the 
	inequalities $\aleph_1=\gma<\gma^+(\aleph_1)=\continuum$ 
	(\Corof{aplus-small-large}) and 
	$\gma=\gma^+(\aleph_1)<\continuum$ (\Thmof{aplus-small-small}) are both 
	consistent.  

	We also give several constructions of mad families with some additional 
	properties. 
\end{abstract}

\section{Introduction}
\institute{Saka\'e Fuchino\at Graduate School of System Informatics, 
  Kobe University, Kobe, Japan.\\
  {\tt fuchino@diamond.kobe-u.ac.jp} 
  \and 
  Stefan Geschke\at Department of Mathematics, University of Hamburg, Germany.\\
  {\tt stefan.geschke@uni-hamburg.de} 
  \and 
  Osvaldo Guzman\at Centre of Mathematics Science, Universidad Nacional 
  Aut\'onoma de M\'exico, Mexico City, Distrito Federal, Mexico.
  {\tt oguzman@matmor.unam.mx}
  \and Lajos Soukup\at Alfr\'ed R\'enyi 
  Institute of Mathematics, Hungarian Academy of Sciences, Budapest, Hungary. {\tt 
	  soukup@renyi.hu}}
Given a family $\calF$ of pairwise almost disjoint countable sets, we can ask 
what the maximal almost disjoint (mad) families 
extending $\calF$ look like. In this note and  
\cite{fuchino-geschke-soukup-2}, we address some instances of this question 
and other related problems. 

Let us begin with the definition of some notions and notation about almost 
disjointness we shall use here. 
Two countable sets $A$, $B$ are said to be {\it almost disjoint\/} ({\em ad} for 
short) if $A\cap B$ is finite. 
A family $\calF$ of countable sets is said to be {\em pairwise almost 
	disjoint\/} ({\em ad\/} for 
short) if any two distinct $A$, $B\in\calF$ are ad. 

If $\calX\subseteq[S]^{\aleph_0}$ and $S=\bigcup\calX$, 
$\calF\subseteq\calX$ is said to be {\em mad in $\calX$} if $\calF$ is ad and 
there is no ad $\calF'$ \st\ $\calF\subsetneqq\calF'\subseteq\calX$. Thus 
an ad family $\calF$ is mad in $\calX$ if and only if there is no 
$X\in\calX$ which is ad\ from every $Y\in\calF$. 
If $\calF$ is mad in $[S]^{\aleph_0}$ for $S=\bigcup\calF$, we say simply 
that $\calF$ is a mad family (on $S$). $S$ as above is called the {\em underlying 
	set\/} of $\calF$. 

Let
\begin{xitemize}
	\xitem[]
	$\gma(\calX)=\min\setof{\cardof{\calF}}{\cardof{\calF}\geq\aleph_0
		\mbox{ and }\calF\mbox{ is mad in }\calX}$. 
\end{xitemize}
Clearly, 
the cardinal invariant $\gma$ known as the almost disjoint number 
(\cite{blass}) can be characterized as:  
\begin{example}
	$\gma=\gma([S]^{\aleph_0})$ for any countable $S$. 
\end{example}

In this paper we concentrate on the case where the underlying set
$S=\bigcup\calX$ (or 
$S=\bigcup\calF$) is countable. In \cite{fuchino-geschke-soukup-2} and 
the forthcoming continuation of this paper, we will 
deal with the cases where $S$ may be also uncountable. 

As the countable $S=\bigcup\calX$, we often use $\omega$ or
$T=\fnsp{\omega>}{2}$ where $T$ is considered as a tree growing downwards. That 
is, for $b$, $b'\in T$, we write $b'\leq_T b$ if $b\subseteq b'$. 
Each $f\in\fnsp{\omega}{2}$ induces the (maximal) branch
\begin{xitemize}
	\xitem[] $B(f)=\setof{f\restr n}{n\in\omega}\subseteq T$ 
\end{xitemize}
in $T$. 

In Section \ref{mad-families}, we consider several cardinal invariants of the form
$\gma(\calX)$ for some $\calX\subseteq[T]^{\aleph_0}$. 

For $\calX\subseteq[S]^{\aleph_0}$ with $S=\bigcup\calX$, let 
\begin{xitemize}
	\xitem[]
	$\calX^\perp=\setof{Y\in[S]^{\aleph_0}}{\forall X\in\calX\ 
		\cardof{X\cap Y}<\aleph_0}$. 
\end{xitemize}
If $Y\in\calX^\perp$ we shall say that $Y$ 
is {\em almost disjoint} (ad) {\em to} $\calX$. 

For an ad family $\calF$, let 
\begin{xitemize}
	\xitem[] $\gma^+(\calF)=\gma(\calF^\perp)$. 
\end{xitemize}
For a cardinal $\kappa$, let 
\begin{xitemize}
	\xitem[] $
	\gma^+(\kappa)=\sup\setof{\gma^+(\calF)}{{}
		\calF\mbox{ is an ad family on }\omega\mbox{ of cardinality }\leq\kappa}
	$.
\end{xitemize}
Clearly, $\gma^+(\omega)=\gma$ and 
$\gma^+(\kappa)\leq\gma^+(\lambda)\leq\continuum$ for 
any $\kappa\leq\lambda\leq\continuum$. 
In Section \ref{mad-over-pad} we give several constructions of ad families  
$\calF$ for which $\calF^\perp$ has some particular property. Using these 
constructions,  we show in Section \ref{aplus} that
$\gma^+(\continuum)=\continuum$ (actually we have
$\gma^+(\bar{\gmo})=\continuum$, see \Thmof{aplus-o}) and the consistency of  
the inequalities 
$\gma=\aleph_1<\gma^+(\aleph_1)=\continuum$ (see \Corof{aplus-small-large}).  
We also show the consistency of 
$\gma^+(\aleph_1)<\continuum$ (\Thmof{aplus-small-small}). 

For notions in the theory of forcing, the reader may consult 
\cite{millennium-book} or  
\cite{kunen-book}. 
We mostly follow the notation and conventions set in \cite{millennium-book} and/or  
\cite{kunen-book}. 
In particular, elements of posets $\poP$ are considered in such a way that stronger 
conditions are smaller. We assume that $\poP$-names are 
constructed just as in \cite{kunen-book} for a \po\ $\poP$ but  
we use alphabets with a tilde below them like
$\utilde{a}$, $\utilde{b}$\vspace{-2pt} etc.\ to denote the $\poP$-names 
corresponding  
to the sets $a$, $b$ etc.\ in the generic extension. 
$V$ denotes the ground model (in which we live). The canonical $\poP$-names of 
elements $a$, $b$ etc.\ of $V$ are denoted by the same symbols with hat like
$\hat{a}$, $\hat{b}$ etc. 
For a \po\ $\poP$ (in $V$) we 
use $V^\poP$ to denote a ``generic'' generic extension $V[G]$ of $V$ by some 
$(V,\poP)$-generic filter $G$.  
Thus $V^\poP\models\ \cdots$ is 
synonymous to $\forces{\poP}{\cdots}$ or
$V\models\forces{\poP}{\cdots}$ and a phrase like: ``Let $W=V^\poP$\,'' 
is to be interpreted as saying: ``Let $W$ be a generic extension of $V$ by 
some/any $(V,\poP)$-generic filter''.

For the notation connected to the set theory of reals see 
\cite{tomek-book} and \cite{blass}. 
By $\continuum$ we denote the size of the continuum 
$2^{\aleph_0}$. $\calM$ and $\calN$ are the ideals of meager sets and null 
sets (e.g.\ over the Cantor space $\fnsp{\omega}{2}$ or the Baire space
$\fnsp{\omega}{\omega}$) respectively.   
For $I=\calM$, $\calN$ etc., $\cov(I)$ and $\non(I)$ 
are {\em covering number} and {\em uniformity} of $I$.

For an infinite cardinal $\kappa$ let $\Cohen{\kappa}=\Fn(\kappa,2)$ or, more 
generally $\Cohen{X}=\Fn(X,2)$ for any set $X$.  
$\Cohen{\kappa}$ is the Cohen 
forcing for adding $\kappa$ many Cohen reals. $\random{\kappa}$ denotes 
the random forcing for adding $\kappa$ many random reals. $\random{\kappa}$ 
is the \po\ consisting of Borel sets of positive measure in
$\fnsp{\kappa}{2}$, which 
corresponds to the homogeneous measure algebra of Maharam 
type $\kappa$. 

For a \po\ $\poP=\pairof{\poP,\leq_\poP}$, $X\subseteq\poP$ and $p\in\poP$, 
let 
\begin{xitemize}
\item[] $X\downarrow p=\setof{q\in X}{q\leq_\poP p}$. 
\end{xitemize}
\section{Mad families and almost disjoint numbers}
\label{mad-families}
One of the advantages of using $T=\fnsp{\omega>}{2}$ as the countable underlying set 
is that we can define some natural subfamilies of $[T]^{\aleph_0}$ such as 
$\calO_T$, $\calA_T$, $\calB_T$ below. 

For $X\subseteq T$, let
\begin{xitemize}
	\xitem[]
	$[X]=\setof{f\in\fnsp{\omega}{2}}{B(f)\subseteq X}$, and  
	\xitem[] 
	$\flr{X}=\setof{f\in\fnsp{\omega}{2}}{\cardof{B(f)\cap X}=\aleph_0}$.
\end{xitemize}
Clearly, we have $[X]\subseteq\flr{X}$. For $X\subseteq T$, let $X^{\uparrow}$ 
be the upward closure of $X$, that is:
\begin{xitemize}
	\xitem[] $X^{\uparrow}=\setof{t\restr n}{t\in X,\,n\leq\ell(t)}$. 
\end{xitemize}
Then we have $\flr{X}\subseteq[X^{\uparrow}]$ for any $X\subseteq T$. 
\begin{definition}[Off-binary sets, \cite{leathrum}] 
	Let
	\begin{xitemize}
	\item[] $\calO_T=\setof{X\in[T]^{\aleph_0}}{
		\flr{X}=\emptyset}$.
	\end{xitemize}
\end{definition}
T.\ Leathrum \cite{leathrum} called elements of $\calO_T$ off-binary sets. Note 
that $\flr{X}=\emptyset$ if and only if there is no branch in $T$ with 
infinite intersection with $X$. 
\begin{definition}[Antichains]
	Let 
	\begin{xitemize}
	\item[] $\calA_T=\setof{X\in[T]^{\aleph_0}}{
		X\mbox{ is an antichain in }T}$.
	\end{xitemize}
\end{definition}
\noindent
Clearly, we have $\calA_T\subseteq\calO_T$. 

Using 
the notation above,  the 
cardinal invariants $\gmo$ and $\bar{\gmo}$ introduced by Leathrum 
\cite{leathrum} can be characterized as:
\begin{xitemize}
	\xitem[] $\gmo=\gma(\calO_T)$, 
	\xitem[] $\bar{\gmo}=\gma(\calA_T)$
\end{xitemize}
(see \cite{leathrum}). 
\iftesting
\mynote{$\gma(\calA_T)=\gma(\calA_{T_\omega})$}
\fi
Leathrum also showed $\gma\leq\gmo\leq\bar{\gmo}$. 
J.\ Brendle \cite{brendle0} proved $\non(\calM)\leq\gmo$. 
\begin{definition}[Sets without infinite antichains]
Let
\begin{xitemize}
\item[] $\calB_T=\setof{X\in [T]^{\aleph_0}}{X
	\mbox{ does not contain any infinite antichain}}$. 
\end{xitemize}
\end{definition}
Note that $\calB_T={\calA_T}^\perp$. 
Elements of $\calB_T$ are those infinite subsets of $T$ which can be 
covered by finitely may branches:  
\begin{lemma}[K.\ Kunen] 
	\label{kunen's lemma}
	Let $X\in [T]^{\aleph_0}$. Then $X\in\calB_T$ 
	if and only if $X$ is covered by finitely may branches in $T$.
\end{lemma}

\begin{proof} If $X$ is covered by finitely many branches in $T$ then $X$ clearly 
does not contain any infinite antichain since otherwise one of the finitely 
many branches would contain an infinite antichain.  

\iffalse
Suppose now that 
$X$ does not contain any infinite antichain. 
Let 
\begin{xitemize}
	\xitem[] $Y=\setof{b\in X}{\mbox{there are incompatible }a,a'\in X
		\mbox{ \st\ }a,a'\leq_T b}$
\end{xitemize}
For each $b\in Y$ let $a^b_0$, $a^b_1\in X$ be \st\ $a^b_0$ and $a^b_1$ are 
incompatible and $a^b_0$, $a^b_1\leq_T b$. 

We claim that $Y$ is finite. 
Otherwise, since $Y$ contains no infinite antichain by
$Y\subseteq X$, there is an infinite branch $C\subseteq Y$ by K\"onig's 
lemma. For each $b\in C$ we can choose $i_b\in 2$ \st\
$a^b_{i_b}\not\in C$. Then $\setof{a^b_{i_b}}{b\in C}$ would be an infinite 
antichain in $X$. But this is a contradiction. 

Let 
\begin{xitemize}
\item[] $D=\setof{d\in X\setminus Y}{d\mbox{ is maximal in }X\setminus Y
	\mbox{ \wrt\ }\leq_T}$. 
\end{xitemize}
Then $D$ is pairwise incompatible and hence finite. 

By definition of $Y$, 
\begin{xitemize}
\item[] $(X\setminus Y)\downarrow d=\setof{b\in X\setminus Y}{b\leq_T d}$ 
\end{xitemize}
is linearly ordered by
$\leq_T$ for each $d\in D$. Thus $X$ is the union of the finitely many linearly 
ordered subsets of $T$: 
\begin{xitemize}
\item[] $(X\setminus Y)\downarrow d,\, d\in D$\ \ and\ \  $\ssetof{b},\, b\in Y$.
\end{xitemize}
This proves the lemma since 
each of these linearly ordered subsets of $T$ can be extended to a branch 
in $T$. 

Suppose for contradiction that $X$ does not contain any infinite 
antichain but neither covered by finitely many branches. 

For $n\in\omega$, let $X_n\in[X]^{\aleph_0}$, $f_n\in\fnsp{\omega}{2}$ and
$t_n\in T$ be defined inductively as follows:
\begin{xitemize}
	\xitem[e-0] $X_0=X$;
	\xitem[e-1] $X_n$ is not a union of finitely many branches;
	\xitem[e-2] $f_n\in\flr{X_n}$;
	\xitem[e-3] $t_n\in X_n$, $t_n\not\in B(f_n)$;
	\xitem[e-4] $X_{n+1}=X_n\downarrow t_n$.
\end{xitemize}
$X_0$ as in \xitemof{e-0} satisfies \xitemof{e-1} by the assumption 
that $X$ is not covered by finitely many branches. \xitemof{e-2} is 
possible by Fodor's lemma and since $X_n\subseteq X$ does not contain any 
infinite antichain. We can find $t_n$ satisfying \xitemof{e-3} by 
\xitemof{e-1}. $t_n$ can be chosen so that $X_{n+1}$ in \xitemof{e-4} 
also satisfies \xitemof{e-1} because the set 
\begin{xitemize}
\item[] $\setof{k\in\omega}{\mbox{there is }u\in X_n\mbox{\ \st\ }
	u\restr k\in B(f_n)\mbox{ but\ }u\restr k+1\not\in B(f_n)}$
\end{xitemize}
is finite: the latter assertion holds since $X_n\subseteq X$ contains no 
infinite antichain. 

Now, for $n\in\omega$, let $m_n\in\omega$ be \st\ $m_n\geq \ell(t_n)$ and
$f_n\restr m_n\in X_n$. Then $\setof{f_n\restr m_n}{n\in\omega}$ is an 
infinite antichain in $X$. This is a contradiction. 
\else
Suppose now that $X$ cannot be covered by finitely many branches. 
By induction on $n$, we choose $t_n\in 2^n$ \st\ 
$t_0=\emptyset$, $t_{n+1}=t_n\concat i$ for some $i\in 2$ and 
\begin{xitemize}
	\xitem[e-5] $X_{n+1}=X\downarrow t_{n+1}$ can not be covered by finitely 
	many branches. 
\end{xitemize}
This is possible 
since $X_0=X$ and $X_n\subseteq (X_n \downarrow (t_n\concat 0)) 
\cup (X_n \downarrow (t_n\concat 1))
\cup \ssetof{t_n}$. 

By \xitemof{e-5}, the branch 
$B=\setof{t_n}{n<\omega}$ does not cover $X_n$ for each $n\in\omega$. So 
we can pick $s_n\in X_n\setminus B$. Let $S=\setof{s_n}{n\in\omega}$.
$S$ is an infinite subset of $X$ since $\ell(s_n)\geq n$ for all $n\in\omega$.
If $C$ is a branch in $T$ different from $B$ then $t_n\notin C$ for some
$n\in\omega$ and so $s_m\notin C$ for all $m\ge n$. Hence $S\cap C$ is finite. 
Moreover  $S\cap B =\emptyset$. So we have 
$\flr{S}=\emptyset$. Thus $S$ should contain an infinite antichain
by K\"onig's Lemma.  
\fi
\qed
\end{proof}
\begin{theorem}[K.\ Kunen]
	\label{Kunen's thm A} 
	$\gma(\calB_T)=\continuum$. 
\end{theorem}
\begin{proof} 
Suppose that $\calF\subseteq\calB_T$ is an ad family of cardinality $<\continuum$. 
We show that $\calF$ is not mad. 
For each $X\in\calF$  
there is $b_X\in[\fnsp{\omega}{2}]^{<\aleph_0}$ \st\
$X\subseteq\bigcup_{f\in b_X}B(f)$ by 
\Lemmaof{kunen's lemma}.  
Since $\calS=\bigcup\setof{b_X}{X\in\calF}$ has cardinality
$\leq\cardof{\calF}\cdot\aleph_0<\continuum$,  
there is $f^*\in\fnsp{\omega}{2}\setminus\calS$. We have $B(f^*)\in\calB_T$ 
and $B(f^*)$ is ad\ to $\calF$. 
\qed
\end{proof}
Let us say $X\subseteq T$ is {\em nowhere dense} if $\flr{X}$ is nowhere 
dense in the Cantor space $\fnsp{\omega}{2}$. It can be easily shown 
that $X$ is nowhere dense if and only if 
\begin{xitemize}
	\xitem[e-6] $\forall t\in T\ \exists t'\leq_T t\ \forall t''\leq_T t'\  
	(t''\not\in X)$. 
\end{xitemize}
Note that, if $X\subseteq T$ is not nowhere dense, then $X$ is dense below 
some $t\in T$ (in terms of forcing). Also note that from \xitemof{e-6} it 
follows that the property of being nowhere dense is absolute. 
\begin{definition}[Nowhere dense sets]
	Let
	\begin{xitemize}
	\item[] 
		$\calND_T=\setof{X\in[T]^{\aleph_0}}{X\mbox{ is nowhere dense\,}}$. 
	\end{xitemize}
\end{definition}

Note that, for $X\in[T]^{\aleph_0}$ with $X=\setof{t_n}{n\in\omega}$, we 
have 
\begin{xitemize}
\item[] $\flr{X}=\bigcap_{n\in\omega}\bigcup_{m>n}[T\downarrow t_m]$. 
\end{xitemize}
In 
particular $\flr{X}$ is a $G_\delta$ subset of $\fnsp{\omega}{2}$. Hence by 
Baire Category Theorem we have 
\begin{xitemize}
\item[]
	$\calND_T=\setof{X\in[T]^{\aleph_0}}{\flr{X}\mbox{ is a meager subset of }
	\fnsp{\omega}{2}}$. 
\end{xitemize}
\begin{lemma}
	\label{claim-0}
	If $X\in[T]^{\aleph_0}$ then there is
	$X'\in[X]^{\aleph_0}$ \st\ $X'\in\calND_T$.
\end{lemma}
\begin{proof} If $\flr{X}=\emptyset$ then $X\in\calND_T$. Thus we can put
$X'=X$. Otherwise let $f\in\flr{X}$ and let $X'=X\cap B(f)$. 
\qed
\end{proof}
\begin{theorem}
	\label{meager-set}
	$\cov(\calM)$, $\gma\leq \gma(\calND_T)$. 
\end{theorem}
\begin{proof} 
For the inequality $\cov(\calM)\leq \gma(\calND_T)$,  
suppose that $\calF\subseteq\calND_T$ is an ad family of cardinality $<\cov(\calM)$.  
Then $\bigcup\setof{\flr{X}}{X\in\calF}\not=\fnsp{\omega}{2}$. Let
$f\in\fnsp{\omega}{2}\setminus\bigcup\setof{\flr{X}}{X\in\calF}$. Then
$B(f)\in\calND_T$ and $B(f)$ is ad\ from all $X\in\calF$. 

To show $\gma\leq \gma(\calND_T)$ suppose that $\calF\subseteq\calND_T$ is an 
ad family of cardinality $<\gma$. Then $\calF$ is not a mad family in
$[T]^{\aleph_0}$. Hence there is some $X\in[T]^{\aleph_0}$ ad\ 
to $\calF$. By \Lemmaof{claim-0}, there is $X'\subseteq X$ \st\
$X'\in\calND_T$. Since $X'$ is also ad\ to $\calF$, it follows 
that $\calF$ is not mad in $\calND_T$. 
\qed
\end{proof}

Let $\sigma$ be the measure on Borel sets of the Cantor space
$\fnsp{\omega}{2}$ defined as the product measure of the probability 
measure on $2$. For $X\subseteq T$, let $\mu(X)=\sigma(\flr{X})$. 
\begin{definition}[Null sets] Let 
	\begin{xitemize}
	\item[] $\calN_T=\setof{X\in[T]^{\aleph_0}}{\mu(X)=0}$. 
	\end{xitemize}
\end{definition}
\begin{theorem}
	\label{null-sets}
	$\cov(\calN)$, $\gma\leq\gma(\calN_T)$.
\end{theorem}
\begin{proof} Similarly to the proof of \Thmof{meager-set}. 
\qed
\end{proof}
\begin{definition}[Nowhere dense null sets]
	Let
	\begin{xitemize}
	\item[] $\calNDN_T=\calND_T\cap\calN_T$. 
	\end{xitemize}
\end{definition}

\begin{lemma}
$\gma(\calND_T)\leq\gma(\calNDN_T)$\ \ and\ \  $\gma(\calN_T)\leq\gma(\calNDN_T)$. 
\end{lemma}
\begin{proof} For the first inequality, suppose that $\calF$ is a mad family in
$\calNDN_T$. Then $\calF$ is an ad family in $\calND_T$. It is also mad in
$\calND_T$. Suppose not. Then there is an $X\in\calND_T$ ad to $\calF$. 
Let $X'\in[X]^{\aleph_0}$ be as in the measure analog of \Lemmaof{claim-0}. Then 
$X'\in\calNDN_T$. Hence $\calF$ is not mad in $\calNDN_T$. This is a 
contradiction. The second inequality can be also proved similarly. 
\qed
\end{proof}
The diagram Fig.\,\ref{fig:1} summarizes the inequalities 
obtained in this section integrated into the cardinal diagram given in Brendle 
\cite{brendle}.
``$\kappa\,\rightarrow\,\lambda$'' in the diagram means that  
``$\kappa\leq\lambda$ is provable in \ZFC''. There are still some open 
questions concerning the 
(in)completeness of this diagram. In particular:
\begin{figure}
\mbox{}\hspace{-11ex}\scalebox{1.35}{\resizebox{0.75\textwidth}{!}{%
\unitlength 1cm
\begin{picture}( 14.9000,  5.7125)( -0.2000, -8.4750)
\put(2.5000,-8.5000){\makebox(0,0){\large$\cov(\calN)$}}%
\put(4.5000,-6.5000){\makebox(0,0){\large$\non(\calM)$}}%
\put(4.5000,-8.5000){\makebox(0,0){\large$\gmb$}}%
\put(8.5000,-6.5000){\makebox(0,0){\large$\gma$}}%
\put(10.0250,-3.5000){\makebox(0,0){\large$\bar{\gmo}$}}%
\put(10.5000,-5.0000){\makebox(0,0){\large$\gmo$}}%
\put(13.0000,-8.5000){\makebox(0,0){\large$\cov(\calM)$}}%
\put(13.0250,-6.5000){\makebox(0,0){\large$\gmd$}}%
\put(15.5250,-6.5000){\makebox(0,0){\large$\non(\calN)$}}%
\put(3.1750,-4.6750){\makebox(0,0)[rb]{\large$\gma(\calN_T)$}}%
%
\special{pn 13}%
\special{pa 1083 3180}%
\special{pa 1693 2707}%
\special{fp}%
\special{sh 1}%
\special{pa 1693 2707}%
\special{pa 1629 2732}%
\special{pa 1652 2740}%
\special{pa 1654 2763}%
\special{pa 1693 2707}%
\special{fp}%
%
\special{pn 13}%
\special{pa 2048 2432}%
\special{pa 2441 2117}%
\special{fp}%
\special{sh 1}%
\special{pa 2441 2117}%
\special{pa 2377 2142}%
\special{pa 2400 2149}%
\special{pa 2402 2173}%
\special{pa 2441 2117}%
\special{fp}%
%
\special{pn 13}%
\special{pa 2727 1890}%
\special{pa 3780 1526}%
\special{fp}%
\special{sh 1}%
\special{pa 3780 1526}%
\special{pa 3711 1529}%
\special{pa 3731 1543}%
\special{pa 3725 1566}%
\special{pa 3780 1526}%
\special{fp}%
\put(6.5250,-4.9750){\makebox(0,0){\large$\gma_\gms$}}%
%
\special{pn 13}%
\special{pa 3219 2550}%
\special{pa 1152 1930}%
\special{fp}%
\special{sh 1}%
\special{pa 1152 1930}%
\special{pa 1209 1967}%
\special{pa 1202 1944}%
\special{pa 1221 1930}%
\special{pa 1152 1930}%
\special{fp}%
\special{pa 975 3189}%
\special{pa 985 1910}%
\special{fp}%
\special{sh 1}%
\special{pa 985 1910}%
\special{pa 964 1976}%
\special{pa 985 1962}%
\special{pa 1003 1976}%
\special{pa 985 1910}%
\special{fp}%
\special{pa 3268 2501}%
\special{pa 2688 2067}%
\special{fp}%
\special{sh 1}%
\special{pa 2688 2067}%
\special{pa 2728 2123}%
\special{pa 2730 2099}%
\special{pa 2751 2091}%
\special{pa 2688 2067}%
\special{fp}%
\special{pa 3426 2501}%
\special{pa 4046 2048}%
\special{fp}%
\special{sh 1}%
\special{pa 4046 2048}%
\special{pa 3981 2070}%
\special{pa 4003 2078}%
\special{pa 4004 2102}%
\special{pa 4046 2048}%
\special{fp}%
%
\special{pn 13}%
\special{pa 2067 2569}%
\special{pa 3957 2018}%
\special{fp}%
\special{sh 1}%
\special{pa 3957 2018}%
\special{pa 3888 2017}%
\special{pa 3907 2033}%
\special{pa 3900 2056}%
\special{pa 3957 2018}%
\special{fp}%
%
\special{pn 13}%
\special{pa 4056 1861}%
\special{pa 3878 1516}%
\special{fp}%
\special{sh 1}%
\special{pa 3878 1516}%
\special{pa 3891 1583}%
\special{pa 3902 1563}%
\special{pa 3926 1565}%
\special{pa 3878 1516}%
\special{fp}%
%
\special{pn 13}%
\special{pa 3544 2530}%
\special{pa 4705 1969}%
\special{fp}%
\special{sh 1}%
\special{pa 4705 1969}%
\special{pa 4637 1980}%
\special{pa 4658 1992}%
\special{pa 4655 2015}%
\special{pa 4705 1969}%
\special{fp}%
%
\special{pn 13}%
\special{pa 4961 3229}%
\special{pa 4154 2077}%
\special{fp}%
\special{sh 1}%
\special{pa 4154 2077}%
\special{pa 4176 2142}%
\special{pa 4185 2121}%
\special{pa 4208 2120}%
\special{pa 4154 2077}%
\special{fp}%
%
\special{pn 13}%
\special{pa 1782 3199}%
\special{pa 1802 2668}%
\special{fp}%
\special{sh 1}%
\special{pa 1802 2668}%
\special{pa 1780 2733}%
\special{pa 1800 2720}%
\special{pa 1819 2734}%
\special{pa 1802 2668}%
\special{fp}%
%
\special{pn 13}%
\special{pa 1841 3268}%
\special{pa 3239 2658}%
\special{fp}%
\special{sh 1}%
\special{pa 3239 2658}%
\special{pa 3171 2666}%
\special{pa 3190 2679}%
\special{pa 3187 2702}%
\special{pa 3239 2658}%
\special{fp}%
%
\special{pn 13}%
\special{pa 1880 3337}%
\special{pa 5020 2609}%
\special{fp}%
\special{sh 1}%
\special{pa 5020 2609}%
\special{pa 4952 2605}%
\special{pa 4969 2621}%
\special{pa 4961 2643}%
\special{pa 5020 2609}%
\special{fp}%
%
\special{pn 13}%
\special{pa 4991 3189}%
\special{pa 4784 1999}%
\special{fp}%
\special{sh 1}%
\special{pa 4784 1999}%
\special{pa 4776 2066}%
\special{pa 4793 2051}%
\special{pa 4814 2060}%
\special{pa 4784 1999}%
\special{fp}%
%
\special{pn 13}%
\special{pa 5109 3189}%
\special{pa 5788 2648}%
\special{fp}%
\special{sh 1}%
\special{pa 5788 2648}%
\special{pa 5724 2674}%
\special{pa 5747 2681}%
\special{pa 5749 2704}%
\special{pa 5788 2648}%
\special{fp}%
%
\special{pn 13}%
\special{pa 5040 3180}%
\special{pa 5099 2658}%
\special{fp}%
\special{sh 1}%
\special{pa 5099 2658}%
\special{pa 5072 2721}%
\special{pa 5093 2710}%
\special{pa 5111 2725}%
\special{pa 5099 2658}%
\special{fp}%
\put(12.9500,-4.8250){\makebox(0,0){\large$\gma(\calND_T)$}}%
%
\special{pn 13}%
\special{pa 4676 1890}%
\special{pa 1516 1290}%
\special{fp}%
\special{sh 1}%
\special{pa 1516 1290}%
\special{pa 1577 1321}%
\special{pa 1567 1300}%
\special{pa 1584 1283}%
\special{pa 1516 1290}%
\special{fp}%
\put(3.7000,-2.9750){\makebox(0,0){\large$\gma(\calNDN_T)$}}%
%
\special{pn 13}%
\special{pa 1054 1625}%
\special{pa 1388 1270}%
\special{fp}%
\special{sh 1}%
\special{pa 1388 1270}%
\special{pa 1329 1305}%
\special{pa 1352 1309}%
\special{pa 1358 1331}%
\special{pa 1388 1270}%
\special{fp}%
\end{picture}%
}}
\caption{}
\label{fig:1}       
\end{figure}
\begin{problem}
	\assert{a} Are the inequalities between $\gma(\calN_T)$, 
	$\gma(\calND_T)$, $\gma(\calNDN_T)$ consistently strict and complete?\smallskip

\assert{b} Are $\gma(\calND_T)$ etc.\ independent from $\gmo$, $\bar{\gmo}$,
$\gma_\gms$ ?
\end{problem}

\section{Ad families $\calF$ for which $\calF^\perp$ is contained in a 
	certain subfamily of $[T]^{\aleph_0}$} 
\label{mad-over-pad}
In this section we give several constructions of ad families with the 
property that the 
sets ad to them in a given generic extension are necessarily in a certain subfamily 
of $[T]^{\aleph_0}$. The constructions 
in this section are used in the proof of some results in the next sections. 

\begin{theorem}
  \label{osvaldo}
  There is an ad family $\calF\subseteq\calA_T$ of cardinality $\non(\calM)$ \st,  
  for any \po\ $\poP$ preserving the non-meagerness of ground-model non-meager sets, we have
  \begin{xitemize}
  \xitem[a-3] $\forces{\poP}{\calF^\perp\subseteq\calND_T}$. 
  \end{xitemize}
\end{theorem}
The following assertion was originally proved under \CH: 
\begin{corollary}
  \label{osvaldo-0}
  There is an ad family $\calF\subseteq\calA_T$ of 
  cardinality $\non(\calM)$ \st, for any cardinal $\kappa$, we have 
	\begin{xitemize}
	\xitem[cc-0] $V^{\Cohen{\kappa}}\models\calF^\perp\subseteq\calND_T$. 
	\end{xitemize}
\end{corollary}
\prf The corollary follows from \Thmof{osvaldo} since the Cohen forcing
$\Cohen{\kappa}$ preserves the non-meagerness of ground-model non-meager sets (see 
e.g.\ 11.3 in \cite{blass})\qed
\qedskip

For the proof of \Thmof{osvaldo}, we use the following lemma. 

Let
\begin{xitemize}
\xitem[] 
  $\calP=\setof{f}{\mapping{f}{X}{\omega}\mbox{ for some }X\in[\omega]^{\aleph_0}}$.
\end{xitemize}
\begin{lemma}
  \label{osvaldo-1}
  There is a mapping $\mapping{F}{\fnsp{\omega}{\omega}}{\calP}$ \st\ 
  \begin{xitemize}
  \xitem[a-4] If $f$, $g\in\fnsp{\omega}{\omega}$, $f\not=g$, then
    $\cardof{F(f)\cap F(g)}<\aleph_0$. 
  \xitem[a-5] If $h\in\fnsp{\omega}{\omega}$ 
    and $X\subseteq\fnsp{\omega}{\omega}$ is non-meager, then there is $f\in X$ 
    \st\ $\cardof{h\cap F(f)}=\aleph_0$. 
  \end{xitemize}
  Furthermore, $F$ as above can be chosen \st\ it is definable and
  absolute in the sense that \xitemof{a-4} and \xitemof{a-5} hold for the 
  extension of $F$ with the  
  same definition in any generic extension of the ground model.
\end{lemma}
\prf Let $\seqof{s_n}{n\in\omega}$ be a one to one recursive enumeration of
$\fnsp{\omega>}{\omega}$. 

For $f\in\fnsp{\omega}{\omega}$, let
$\dom(F(f))=\setof{n\in\omega}{s_n\subseteq f}$. Let  
$\mapping{F(f)}{\dom(F(f))}{\omega}$ be defined by 
\begin{xitemize}
\xitem[a-6] $F(f)(n)=f(\cardof{s_n})$
\end{xitemize}
for $n\in\dom(F(f))$. 
\begin{claim}
  This $F$ is as desired. 
\end{claim}
\prfofClaim
It is clear that $F$ satisfies \xitemof{a-4} --- note that it is crucial here 
that the enumeration $\seqof{s_n}{n\in\omega}$ is chosen to be one to one. 

To show that $F$ also satisfies \xitemof{a-5}, 
suppose $h\in\fnsp{\omega}{\omega}$. It is enough to show that 
\begin{xitemize}
\xitem[a-7] 
  $N(h)=\setof{g\in\fnsp{\omega}{\omega}}{\cardof{h\cap F(g)}<\aleph_0}$ is 
  a meager subset of $\fnsp{\omega}{\omega}$. 
\end{xitemize}
For $k\in\omega$, let
$N_k(h)=\setof{g\in\fnsp{\omega}{\omega}}{\cardof{h\cap F(g)}<k}$. 

Since $N(h)=\bigcup_{k\in\omega}N_k(h)$, it is enough to show that 
$N_k(h)$ is a nowhere dense subset of $\fnsp{\omega}{\omega}$ for each $k\in\omega$. 

For this, we prove, by induction on $k$,  
\begin{xitemize}
\xitem[a-8] 
  For any $s\in\fnsp{\omega>}{\omega}$, there are $s'\in\fnsp{\omega>}{\omega}$ 
  and $m'\in\omega$ \st\ 
  \st\ $s'\subseteq s$ and $\cardof{(h\restr m')\cap F(g)}\geq k$ for all $g\in[s']$. 
\end{xitemize}
Suppose that \xitemof{a-8} holds for $k=\ell$ and let $s\in\fnsp{\omega>}{\omega}$. 
By the induction hypothesis we may assume \wolog\ that there is an $m\in\omega$ 
\st\ 
$\cardof{(h\restr m)\cap F(f)}\geq\ell$ for all $g\in[s]$. 

Let $n\in\omega$ be \st\ $n\geq m, \cardof{s}$ and $s_n\supseteq s$. Let
\begin{xitemize}
\xitem[] 
  $s'=s_n\cup\ssetof{\pairof{\cardof{s_n},h(n)}}$. 
\end{xitemize}
For any $g\in[s']$, we have $n\in\dom(F(g))$ by $s_n\subseteq s'\subseteq g$, and
$F(g)(n)=g(\cardof{s_n})=h(n)$. Letting $m'=n+1$, we have  
$\cardof{(h\restr m')\cap F(g)}\geq\ell+1$. 
Thus, \xitemof{a-8} holds for $k=\ell+1$ with these $s'$ and $m'$. 
\qedofClaim
\qedskip

The definability and the absoluteness of $F$ is clear from the 
construction given above. 
\qed
\qedskip

\noindent
{\bf Proof of \bfThmof{osvaldo}:}\ \  Let
\begin{xitemize}
\xitem[a-9] 
  $Q=\setof{q\in T}{q(n)\mbox{ is eventually }0}$. 
\end{xitemize}
That is, for
$q\in T$, $q\in Q$ if and only if $\cardof{\setof{n\in\omega}{q(n)=1}}<\aleph_0$. 

For $q\in Q$, let 
\begin{xitemize}
\xitem[a-10] 
  $\ell_q=\min\setof{\ell\in\omega}{
  \forall m\ (\ell\leq m\,\rightarrow\,q(m)=0)}$.
\end{xitemize}
Let $\seqof{q_n}{n\in\omega}$ be a one to one enumeration of $Q$. 

For $n$, $k\in\omega$ let 
\begin{xitemize}
\xitem[a-11] 
  $T_{n,k}=\setof{s\in T}{q_n\restr(\ell_q+k)\cup\ssetof{\pairof{\ell_q+k,1}}
  \subseteq s}$
\end{xitemize}
and let $\seqof{s_{n,k,i}}{i\in\omega}$ be a one to one enumeration of $T_{n,k}$. 
Let $F$ be as in \Lemmaof{osvaldo-1}. For $n\in\omega$ and
$f\in\fnsp{\omega}{\omega}$, let 
\begin{xitemize}
\xitem[a-12] 
  $F_n(f)=\setof{s_{n,k,i}}{k\in\dom(F(f)),\,i=F(f)(k)}$.
\end{xitemize}

Let $N\subseteq\fnsp{\omega}{\omega}$ be a non-meager set with
$\cardof{N}=\non(\calM)$. Let $\calF_n=F_n\imageof N$ and
$\calF=\bigcup_{n\in\omega}\calF_n$. 

We show that this $\calF$ is as desired:
\begin{claim}
  \assertof{1} $\calF\subseteq\calA_T$.\smallskip

  \assertof{2} $\calF$ is ad.\smallskip

  \assertof{3} \xitemof{a-3} holds for all \po\ $\poP$ preserving non-meagerness 
  of ground-model non-meager sets. 
\end{claim}
\prfofClaim
\assertof{1}: Suppose that $A\in\calF$ and $A=F_n(f)$ for some $n\in\omega$ and
$f\in N$. If $s_0$, $s_1$ are two different elements of $A$, then there 
are $k_0$, $k_1\in\dom(F(f))$, $k_0\not=k_1$ and $i_0$, $i_1\in\omega$ \st\
$s_0=s_{n,k_0,i_0}$ and $s_1=s_{n,k_1,i_1}$. Since $s_0\in T_{n,k_0}$ and
$s_1\in T_{n,k_1}$, it follows that $s_0$ and $s_1$ are incompatible. \smallskip

\assertof{2}: Suppose that $A_0$, $A_1\in\calF$ with $A_0\not=A_1$. Let 
$A_0=F_{n_0}(f_0)$ and $A_1=F_{n_1}(f_1)$. If $n_0\not=n_1$ then we have
$\cardof{A_0\cap A_1}\leq 1$. Then $f_0\not=f_1$. Thus, by \xitemof{a-4}, 
$\cardof{A_0\cap A_1}=\cardof{F(f_0)\cap F(f_1)}<\aleph_0$. \smallskip

\assertof{3}: Let $G$ be a $(V,\poP)$-generic set and we work in $V[G]$. Note, 
that by our assumption, $N$ is still non-meager in $V[G]$. 

Suppose that $B\in [T]^{\aleph_0}\setminus\calND_T$. We have to show that
$\cardof{A\cap B}=\aleph_0$ for some $A\in\calF$. 

Since $B\not\in\calND_T$ there 
is $n\in\omega$ \st\ $B\downarrow (q_n\restr\ell_{q_n})$ is dense below
$q_n\restr\ell_{q_n}$. It follows that, for each $k\in\omega$, there is 
$h(k)\in\omega$ \st\ $s_{n,k,h(k)}\in B$. By \xitemof{a-5} (which still holds in 
the generic extension $V[G]$), there is $f\in M$ \st\
$\cardof{h\cap F(f)}=\aleph_0$. 

By the definition of $h$ and $F_n(f)$, it follows that
$\cardof{B\cap F_n(f)}=\aleph_0$. 
\qedofClaim
\qed
\qedskip

We can also obtain a variation of \Thmof{osvaldo} if 
our ground 
model is a generic extension of some inner model by adding uncountably may 
Cohen reals. Note that $\non(\calM)=\aleph_1$ holds in such a ground model.

\begin{theorem}
	\label{cohen-nd-1}
	Suppose that $W=V^\Cohen{\omega_1}$. Then, in $W$, there is an ad family
	$\calF\subseteq\calA_T$ of cardinality $\aleph_1$ \st 
  \begin{xitemize}
    \xitem[cohen-nd-a] for any c.c.c.\ \po\ $\poP$ with $\poP\in V$, we have
    $W^\poP\models\calF^\perp\subseteq\calND_T$.  
  \end{xitemize}
\end{theorem}
\begin{proof} Let $A\in[T]^{\aleph_0}\cap V$ be an antichain and let 
$\seqof{t^*_n}{n\in\omega}$ be a one to one enumeration of $A$. 

Let $G$ be a $(V,\Cohen{\omega_1})$-generic filter and $W=V[G]$. 
For $p\in\Cohen{\omega_1}$, $\alpha<\omega_1$ and $k\in\omega$, let  
\begin{xitemize}
\item[] $f^p_\alpha=\setof{\pairof{n,i}\in\omega\times\omega}{\pairof{\omega\alpha+3n,i}\in p}$;\medskip
\item[] $n^p_{\alpha,k}=\left\{\,
  \begin{array}{@{}ll}
    n, &\mbox{if }[\omega\alpha,\omega\alpha+3n+1]\subseteq\dom(p),\\
    &\phantom{\mbox{if }}
    p(\omega\alpha+3n+1)=1\mbox{ and}\\
    &\phantom{\mbox{if }}
    \cardof{\setof{m<n}{p(\omega\alpha+3m+1)=1}}=k,\\[2\jot]
    \mbox{undefined}, &\mbox{if there is no such }n\mbox{ as above;}
  \end{array}\right.
  $\medskip
\item[] $t^p_\alpha=\left\{\,
  \begin{array}{@{}ll}
    \setof{\pairof{n,i}\in\omega\times\omega}{
      n<n^p_{\alpha,0},\,\pairof{\omega\alpha+3n+2,i}\in p},\\[\jot]
    \phantom{\mbox{undefined},\qquad}\mbox{if }n^p_{\alpha,0}\mbox{ is defined,}\\[2\jot]
    \mbox{undefined},\qquad\mbox{otherwise}
  \end{array}\right.
  $

  and 
\item[] $t^p_{\alpha,k}=\left\{\,
  \begin{array}{@{}ll}
    \setof{\pairof{n,i}\in\omega\times\omega}{
      n<n^p_{\alpha,k+1},\,\pairof{\omega\alpha+3n+2,i}\in p},\\[\jot]
    \phantom{\mbox{undefined},\qquad}\mbox{if }n^p_{\alpha,k+1}\mbox{ is defined,}\\[2\jot]
    \mbox{undefined},\qquad\mbox{otherwise.}
  \end{array}\right.
  $
\end{xitemize}
Let 
\begin{xitemize}
\item[] $f^G_\alpha=\bigcup_{p\in G}f^p_\alpha$,
\item[] $t^G_{\alpha}=t^p_{\alpha}$ for some $p\in G$ \st\ $t^p_{\alpha}$ 
  is defined, and 
\item[] $t^G_{\alpha,k}=t^p_{\alpha,k}$ for some $p\in G$ \st\ $t^p_{\alpha,k}$ 
  is defined.
\end{xitemize}

For $\alpha\in\omega_1$, let
\begin{xitemize}
  \xitem[A-alpha]
  $A_\alpha=\setof{t^G_\alpha\concat t^*_k\concat t^G_{\alpha,k}}{k\in\omega}$.  
\end{xitemize}

Clearly each $A_\alpha$ is an antichain in $T$. 

$A_\alpha$, $\alpha<\omega_1$ are 
pairwise almost disjoint: Suppose that $\alpha<\beta<\omega_1$. Then there is
$k_0<\omega$ \st\  
$t^G_{\alpha,k}\not=f^G_{\beta,k}$ for all $k\in\omega\setminus k_0$. It follows 
that
$A_\alpha\cap A_\beta\subseteq\setof{t^G_\alpha\concat t^*_k\concat t^G_{\alpha,k}}{k<k_0}$. 

We show that $\calF=\setof{A_\alpha}{\alpha<\omega_1}$ satisfies 
\xitemof{cohen-nd-a}. 

Suppose that $\poP$ is a c.c.c.\ poset (in $W$) and $\poP\in V$. Let $H$ be a
$(W,\poP)$-generic filter. It is enough to show that, in $W[H]$, if 
$X\in[T]^{\aleph_0}$ is not nowhere dense then $X$ is not almost ad to $\calF$. 

By the c.c.c.\ of $\Cohen{\omega_1}*\hat{\poP} \sim \Cohen{\omega_1}\times\poP$, 
there is an $\alpha^*\in\omega_1$ \st\
$X\in V[G\restr\Cohen{\omega\alpha^*}][H]$. Let $t\in T$ be \st\ $X$ is dense 
below $t$. Then 
\[D=\setof{p\in\Cohen{\omega_1\setminus\omega\alpha^*}}{t^p_\alpha\supseteq t
  \mbox{ for some }\alpha\in\omega_1\setminus\omega\alpha^*}
\]\noindent
is dense in $\Cohen{\omega_1\setminus\omega\alpha^*}$. 

For $p\in D$ and $\alpha\in\omega_1\setminus\omega\alpha^*$ \st\
$t^p_\alpha\supseteq t$, letting $\utilde{A}_\alpha$ a
$\Cohen{\omega_1\setminus\alpha^*}$-name of $A_\alpha$, we have 
$p\forces{\Cohen{\omega_1\setminus\omega\alpha^*}}{
  \cardof{\utilde{A}_\alpha\cap X\downarrow t}=\aleph_0}$ 
by \xitemof{A-alpha} and since $X$ is dense below $t$. 

By genericity, it follows that, in $W[G]$, there is $\alpha<\omega_1$ \st\
$\cardof{A_\alpha\cap X}=\aleph_0$. \mbox{}\qed
\end{proof}
A measure version of \Thmabove\ also holds:
\begin{theorem}
	\label{random-n}
	Let $W=V^\Cohen{\omega_1}$. Then, in $W$, there is an ad family $\calF$ 
	in $\calN_T$ of cardinality $\aleph_1$ \st\, for any c.c.c.\ \po\  
	$\poP$ with $\poP\in V$, we have 
	$W^{\poP}\models\calF^\perp\subseteq\calO_T$. 
\end{theorem}
For the proof of \Thmabove\ we note first the following:
\begin{lemma}
	\label{null-set}
	Suppose that $X\subseteq T$ is \st\ $X=\setof{t_k}{k\in\omega}$ 
	for some enumeration $t_k$, $k\in\omega$ of $X$ with
	$\ell(t_k)\geq k$  for all $k\in\omega$. Then $X\in\calN_T$. 
\end{lemma}
\begin{proof} For all $n\in\omega$, we have
$\flr{X}\subseteq\bigcup_{k\in\omega\setminus n}\flr{T\downarrow t_k}$. 
Hence 
\begin{xitemize}
\item[] 
	$\mu(X)
	=\sigma(\flr{X})\leq\sum_{k\in\omega\setminus n}\sigma(\flr{T\downarrow t_k})
	\leq\sum_{k\in\omega\setminus n}2^k=2^{-n}$. 
\end{xitemize}
It follows that $\mu(X)=0$. 
\qed
\end{proof}
\begin{proof}[of \itThmof{random-n}]
Let $G$ be a $(V,\Cohen{\omega_1})$-generic filter and $W=V[G]$. In $W$, let
\begin{xitemize}
	\item[] $f^G_\alpha=\setof{\pairof{n,i}}{\pairof{\omega\alpha+n,i}\in p
		\mbox{ for some }p\in G}$
\end{xitemize}
for $\alpha<\omega_1$ and let
$g^G_\alpha\in\fnsp{\omega}{\omega}$ be the increasing enumeration of
$\left(f^G_\alpha\right)^{-1}[\ssetof{1}]$. 

Further in $W$, we  construct inductively
$A_\alpha\in\calN_T$,  $\alpha<\omega_1$ as follows. 

For $n\in\omega$, let $A_n\in\calN_T$ be \st\ $\seqof{A_n}{n\in\omega}$ is 
a partition of $T$ in $V$. We can be easily find such $A_n$'s by \Lemmaof{null-set}. 

For $\omega\leq\alpha<\omega_1$, suppose that pairwise almost disjoint $A_\beta$,
$\beta<\alpha$ have been constructed. Let $\seqof{B_\ell}{\ell\in\omega}$ 
be an enumeration of $\setof{A_\beta}{\beta<\alpha}$ and, for each
$n\in\omega$, let $\seqof{b_{n,m}}{m\in\omega}$ be an enumeration of
\begin{xitemize}
	\xitem[d-0] $C_n=T\setminus\left(\fnsp{n>}{2}\cup\setof{B_\ell}{\ell<n}\right)$. 
\end{xitemize}
Let 
\begin{xitemize}
	\xitem[d-1] 
	$A_\alpha=\setof{b_{n,g^G_\alpha(n)}}{n\in\omega}$. 
\end{xitemize}
$A_\alpha\in\calN_T$ by \xitemof{d-0} and \Lemmaof{null-set}. $A_\alpha$ is ad to
$\setof{A_\beta}{\beta<\alpha}$ by 
\xitemof{d-0} and \xitemof{d-1}.

We show that $\calF=\setof{A_\alpha}{\alpha<\omega_1}$ is as desired. 
Suppose that $\poP$ is c.c.c.\ (in $W$) and $\poP\in V$. 
Let $H$ be a $(W,\poP)$-generic filter. It is enough to show that, in $W[H]$, if 
$X\in[T]^{\aleph_0}\setminus\calO_T$ then $X$ is not ad to $\calF$. 
So suppose that (in $W[H]$) $X\in [T]^{\aleph_0}\setminus\calO_T$ and
$f\in\flr{X}$. Let $B=X\cap B(f)$. 
By the c.c.c.\ 
of $\Cohen{\omega_1}\ast\hat{\poP}\sim\Cohen{\omega_1}\times\poP$,  there is an  
$\alpha^*\in\omega_1\setminus\omega$ \st\
$B\in V[(G\restr\Cohen{\omega\alpha^*})][H]$.  
If $B\cap A_\alpha$ is infinite for some $\alpha<\alpha^*$ then we are done. So 
assume that $B$ is ad to all $A_\alpha$, $\alpha<\alpha^*$. Then
$B\cap C_n$ is infinite for all $n\in\omega$. 
Since $f^G_{\alpha^*}$ is a 
Cohen real generic over $V[(G\restr\Cohen{\omega\alpha^*})][H]$, 
it follows that $B\cap A_{\alpha^*}$ is infinite. 
\qed
\end{proof}

\section{Almost disjoint numbers over ad families}
In this section we turn to questions on the possible values of $\gma^+(\cdot)$. 
\label{aplus}
\begin{theorem}{\rm (K.\ Kunen)}
	\label{aplus-o}
	$\gma^+(\bar{\gmo})=\continuum$. 
\end{theorem}
\begin{proof}
Let $\calF$ be any mad family in $\calA_T$ of cardinality $\bar{\gmo}$. By 
maximality of $\calF$ we have $\calF^\perp=\calB_T$. If 
$\calG\subseteq[T]^{\aleph_0}$ is disjoint from $\calF$ 
and $\calF\cup\calG$ is mad then $\calG$ is  
mad in $\calB_T$ and hence $\cardof{\calG}=\continuum$ by \Thmof{Kunen's thm A}. 
\qed
\end{proof}
\begin{theorem}
	\label{aplus-small-Cohen}
	$V^{\Cohen{\kappa}}\models\gma^+(\aleph_1)\geq\kappa$ for all regular
	$\kappa$. 
\end{theorem}
\begin{proof} If $\kappa=\omega_1$ this is trivial. So suppose 
that $\kappa>\omega_1$. Let $W=V^{\Cohen{\omega_1}}$. Then
$V^{\Cohen{\kappa}}=W^{\Cohen{\kappa\setminus\omega_1}}$. Let $\calF$ be as 
in the proof of \Thmof{cohen-nd-1}. Suppose that
$\tilde{\calF}\supseteq\calF$ is mad on $T$ in $V^{\Cohen{\kappa}}$. Then
$\tilde{\calF}\subseteq\left(\calND_T\right)^{V^{\Cohen{\kappa}}}$. 
Since $V^\Cohen{\kappa}\models \cov(\calM)\geq\kappa$, it follows that
$\cardof{\tilde{\calF}}\geq\kappa$ by \Thmof{meager-set}. 
\qed
\end{proof}
\begin{corollary}
	\label{aplus-small-large}
	The inequality $\gma=\aleph_1<\gma^+(\aleph_1)=\continuum$ is consistent. 
\end{corollary}
\begin{proof}
Start from a model $V$ of \CH. Since there is 
a $\Cohen{\kappa}$-indestructible mad family in $V$ 
it follows that
$V^{\Cohen{\omega_2}}\models\gma=\aleph_1$ (see e.g.\ 
\cite{kunen-book}, Theorem 2.3).  On the other hand we have 
$V^{\Cohen{\omega_2}}\models\gma^+(\aleph_1)=\aleph_2=\continuum$ by 
	\Thmof{aplus-small-Cohen}. 
\qed
\end{proof}
\begin{theorem}
	\label{aplus-small-small}
	The inequality $\gma^+(\aleph_1)<\continuum$ is consistent. 
\end{theorem}
\noindent
For the proof of the theorem we use the following forcing notions:
for a family
$\calI\subseteq\setof{A\in[\omega]^{\aleph_0}}{\cardof{\omega\setminus A}=\aleph_0}$ 
closed under union, let 
$\poQ_\calI=\pairof{\poQ_\calI,\leq_{\poQ_\calI}}$ be the \po\ defined by 
\begin{xitemize}
\item[] $\poQ_\calI=\Cohen{\omega}\times\calI$\,;
\end{xitemize}
For all $\pairof{s,A}$, $\pairof{s',A'}\in\poQ_\calI$
\begin{xitemize}
	\xitem[] $
	\begin{array}[t]{r@{}l}
		\pairof{s',A'}\leq_{\poQ_\calI}\pairof{s,A}\ \ \Leftrightarrow\ \ 
		&s\subseteq s',\ A\subseteq A'\mbox{ and }\\
		&\mbox{}\hspace{-12pt}\forall n\in\dom(s')\setminus\dom(s)\ (n\in A\ \rightarrow\ s'(n)=0).
	\end{array}
	$
\end{xitemize}
Clearly $\poQ_\calI$ is $\sigma$-centered. 

For a $(V,\poQ_\calI)$-generic $G$, let
\begin{xitemize}
\item[] $f_G=\bigcup\setof{s}{\pairof{s,A}\in G\mbox{ for some }A\in\calI}$ 
	and 
\item[] $A_G=f^{-1}_G\imageof\ssetof{1}$.
\end{xitemize}
Let $\tilde{\calI}$ be the ideal in $[\omega]^{\aleph_0}$ generated from
$\calI$ (i.e.\ the downward closure of $\calI$ \wrt\ $\subseteq$). 
By the genericity of $G$ and the definition of $\leq_{\poQ_\calI}$ it is 
easy to see that $A_G$ is infinite and 
\begin{xitemize}
	\xitem[c-10] for every $B\in([\omega]^{\aleph_0})^V$, $A_G$ is almost disjoint 
	from $B$\ \ $\Leftrightarrow$\ \ 	$B\in\tilde{\calI}$. 
\end{xitemize}
\noindent
\begin{proof}[of \itThmof{aplus-small-small}]
Working in a ground model $V$ of  
$2^{\aleph_0}=2^{\aleph_1}=\aleph_3$, let 
\begin{xitemize}
\item[] 
	$\seqof{\poP_\alpha,\utilde{\poQ}_\beta}{\alpha\leq\omega_2,\,\beta<\omega_2}$
\end{xitemize}
be the finite support iteration of c.c.c.\ \pos\ defined as follows:
for $\beta<\omega_2$, let $\utilde{\poQ}_\beta$ be the $\poP_\beta$-name of 
the finite support (side-by-side) product of 
\begin{xitemize}
	\xitem[c-11] $\poQ_{\tilde{\calF}}$, $\tilde{\calF}\in\Phi$ 
\end{xitemize}
where
\begin{xitemize}
\item[] $
	\begin{array}[t]{r@{}l}
		\Phi=\setof{\tilde{\calF}}{{}&\tilde{\calF}
			\mbox{ is an ideal in }[\omega]^{\aleph_0}\\
			&
			\mbox{ generated from an ad 
				family in }
					 [\omega]^{\aleph_0}\mbox{ of cardinality }\aleph_1}
	\end{array}
	$ 
\end{xitemize}
in $V^{\poP_\beta}$. We have 
\begin{xitemize}
\item[] $V^{\poP_\beta}\models\utilde{\poQ}_\beta\mbox{ satisfies the c.c.c.}$
\end{xitemize}
since 
$V^{\poP_\beta}\models\poQ_{\tilde{\calF}}\mbox{ is }\sigma
	\mbox{-centered for all }\tilde{\calF}\in\Phi$. 
By induction on $\alpha\leq\omega_2$, we can show that $\poP_\alpha$ satisfies the 
c.c.c.\ and $\cardof{\poP_\alpha}\leq 2^{\aleph_1}=\aleph_3$ for all
$\alpha\leq\omega_2$. It follows that 
\begin{xitemize}
	\xitem[] $V^{\poP_{\omega_2}}\models 2^{\aleph_0}=2^{\aleph_1}=\aleph_3$.
\end{xitemize}
Thus the following claim finishes the proof:
\begin{claim}
	$V^{\poP_{\omega_2}}\models\gma=\gma^+(\aleph_1)=\aleph_2$.
\end{claim}
\prfofClaim
Working  in $V^{\poP_{\omega_2}}$, suppose that $\calF$ is an ad family in 
$[\omega]^{\aleph_0}$ of cardinality $\aleph_1$. By the c.c.c.\ of
$\poP_{\omega_2}$, there is some $\alpha^*<\omega_2$ \st\
$\calF\in V^{\poP_{\alpha^*}}$. By \xitemof{c-11} and \xitemof{c-10}, there 
are $A_\alpha$, $\alpha\in\omega_2\setminus\alpha^*$ \st\ 
\begin{xitemize}
	\xitem[] for every $B\in([\omega]^{\aleph_0})^{V^{\poP_{\alpha}}}$, 
$A_\alpha$ is ad from $B$\ \ $\Leftrightarrow$\ \ $B\in$  the ideal 
	generated from $\calF\cup\setof{A_\beta}{\beta\in\alpha\setminus\alpha^*}$.
\end{xitemize}
Since 
$([\omega]^{\aleph_0})^{V^{\poP_{\omega_2}}}
	=\bigcup_{\alpha<\omega_2}([\omega]^{\aleph_0})^{V^{\poP_\alpha}}$, 
it follows that
$\calF\cup\setof{A_\alpha}{\alpha\in\omega_2\setminus\alpha^*}$ is a mad 
family in $V^{\poP_{\omega_2}}$. This shows that
$V^{\poP_{\omega_2}}\models\gma^+(\aleph_1)\leq\aleph_2$. 

We also have $V^{\poP_{\omega_2}}\models\gma\geq\aleph_2$: 
for any ad family $\calG\subseteq([\omega]^{\aleph_0})^{V^{\poP_{\omega_2}}}$ 
of cardinality $\leq\aleph_1$, there is some $\alpha^*<\omega_2$  
\st\ $\calG\in V^{\poP_{\alpha^*}}$. 
But $\smash{\utilde{\poQ}_{\alpha^*}}$\smallskip\ 
adds an infinite subset of $\omega$ almost disjoint to every element of $\calG$. 
Hence $\calG$ is not mad. 
\qedofClaim\qed
\end{proof}
Clearly, the method of the proof of \Thmabove\ cannot produce a model of
$\gma^+(\aleph_1)=\aleph_1<\continuum$.
\begin{problem}
	Is $\gma^+(\aleph_1)=\aleph_1<\continuum$ consistent? 
\end{problem}

All infinite cardinals less than or equal to the continuum $\continuum$ can 
be represented as $\gma^+(\calF)$ for some $\calF$. 
\begin{theorem}
	\label{aplus-all}
	For any infinite $\kappa\leq\continuum$, there is an ad family 
	$\calF\subseteq[T]^{\aleph_0}$ of cardinality $\continuum$ \st\
	$\gma^+(\calF)=\kappa$.  
\end{theorem}
\begin{proof} Let $\calF'$ be a mad family in $\calA_T$. Then by 
\Lemmaof{kunen's lemma}, we have 
\begin{xitemize}
	\xitem[d-2] $\calF'^\perp=\calB_T$. 
\end{xitemize}
Let $X$ and $X'$ be disjoint with
$\fnsp{\omega}{2}=X\cup X'$, $\cardof{X}=\continuum$ and
$\cardof{X'}=\kappa$. Let 
\begin{xitemize}
\item[] $\calF=\calF'\cup\setof{B(f)}{f\in X}$. 
\end{xitemize}
Clearly $\calF$ is an ad family. By \xitemof{d-2} we have
$\calF^\perp\subseteq\calB_T$. 

We claim $\gma^+(\calF)=\kappa$: 
Since $\calF\cup\setof{B(f)}{f\in X'}$ is a mad family by 
\Lemmaof{kunen's lemma}, we have $\gma^+(\calF)\leq\kappa$. Again by 
\Lemmaof{kunen's lemma}, if $\calG\subseteq\calF^\perp$ is an ad family of 
cardinality $<\kappa$, then there is $f\in X'$ \st\ $B(f)$ is ad from 
every $B\in\calG$. Thus $\gma^+(\calF)\geq\kappa$. 
\qed
\end{proof}

\section{Destructibility of mad families}
For a \po\ $\poP$, a mad family $\calF$ in $[T]^{\aleph_0}$ is 
said to be {\em $\poP$-destructible} if 
\begin{xitemize}
	\item[] $V^\poP\models\calF$ is not mad in
$[T]^{\aleph_0}$. 
\end{xitemize}
Otherwise it is {\em$\poP$-indestructible}.  

The results in Section \ref{mad-over-pad} can be also formulated in terms of 
destructibility of mad families. 
\begin{theorem}
	\label{abs-0}
	\assert{1}
	There is an ad family $\calF\subseteq\calA_T$ of size $\non(\calM)$ which cannot be 
	extended to a $\Cohen{\omega}$-indestructible mad family in any generic 
	extension of the ground model $V^{\poP}$ as long as non-meager sets in $V$ 
  remain non-meager in $V^\poP$. \smallskip

	\assert{2}
	Let $W=V^{\Cohen{\omega_1}}$. Then, in $W$, there is an ad family 
	$\calF\subseteq\calND_T$ of cardinality $\aleph_1$ \st, in any generic 
	extension of $W$ by a c.c.c.\ \po\ $\poP$ with $\poP\in V$, $\calF$ 
	cannot be extended to a $\Cohen{\omega}$-indestructible mad family. \smallskip

	\assert{3}
	Let $W=V^{\Cohen{\omega_1}}$. Then, in $W$, there is an ad family 
	$\calF\subseteq\calN_T$ of cardinality $\aleph_1$ \st, in any generic 
	extension of $W$ by a c.c.c.\ \po\ $\poP$ with $\poP\in V$, $\calF$ 
	cannot be extended to a $\random{\omega}$-indestructible mad family. 
\end{theorem}
\begin{proof} \assertof{1}:
The family $\calF$ as in \Thmof{osvaldo} will do. Since we have
$\calF'\subseteq \calND_T$ for any mad
$\calF'$ extending $\calF$ in $V^\poP$, a further Cohen 
real over 
$V^\poP$ introduces a branch almost avoiding all elements of
$\calF'$. Thus $\calF'$ is no longer mad in
$V^{\poP\ast\Cohen{\omega}}$. 
\smallskip

\assertof{2}:
By \Thmof{cohen-nd-1} and by an argument similar to the proof of 
\assertof{1}. 
\smallskip

\assertof{3}: 
In $W$, let $\calF$ be as in the proof of \Thmof{random-n}. Then any 
mad $\calF'\supseteq\calF$ on $T$ in any $W^\poP$ for $\poP$ as above 
is included in $\calN_T$ by
$\calO_T\subseteq\calN_T$. Hence, in $W^{\poP\ast\random{\omega}}$, the 
random real $f$ over $W^\poP$ introduces the branch $B(f)$ almost avoiding 
all elements  
of $\calF'$. Thus $\calF'$ is no longer mad in $W^{\poP\ast\random{\omega}}$. 
\qed
\end{proof}

\section{$\kappa$-almost decided and $\lambda$-minimal mad families}
In this final section we collect several other constructions of mad families 
with some additional properties.

\newcommand{\ide}{\operatorname{\cal I}}
Given an ad family $\calF$ on $T$ let $\ide (\calF)$ be the ideal on
$T$ generated by $\calF\cup [T]^{<{\omega}}$, i.e. for $S\subset T$ we have
$S\in \ide(\calF)$ if $S\subset^*\cup\calF'$ for some  finite subfamily 
$\calF'$ of $\calF$.

Let  $\calF$  be a mad family on $T$ and  $\calB\subseteq\calF$.
Clearly $\calB^\perp \supseteq \ide(\calF\setminus \calB)\setminus[T]^{<\aleph_0}$. 
We say  that $\calB$ {\em almost decides} $\calF$
if $\calB^\perp = \ide(\calF\setminus \calB)\setminus[T]^{<\aleph_0}$.
 A mad family $\calF$ is said to be {\em$\kappa$-almost decided\/} if every 
$\calB\in[\calF]^{\kappa}$ almost decides $\calF$.

\begin{theorem}
	\label{c-almost decided}
	Assume that $\MA(\sigma\mbox{-centered\/{}})$ holds. 
        Then there is a $\continuum$-almost decided mad 
	family $\calF$ on $T$. 
\end{theorem}
\begin{proof} Let $\seqof{B_\beta}{\beta<\continuum}$ be an enumeration of
$[T]^{\aleph_0}$. We define $A_\alpha$, $\alpha<\continuum$  
inductively \st\
\begin{xitemize}
	\xitem[f-0] $\setof{A_n}{n\in\omega}$ is a partition of $T$ into infinite 
	subsets;
\end{xitemize}
For all $\alpha\in\continuum\setminus\omega$
\begin{xitemize}
	\xitem[f-1] $A_\alpha$ is ad from $A_\beta$ for all $\beta<\alpha$;
	\xitem[f-2] For $\beta<\alpha$, if
	$B_{\beta}\notin \ide(\setof{A_\delta}{{\delta}<{\alpha}})$ 
	then $\cardof{A_\alpha\cap B_\beta}=\aleph_0$;
\end{xitemize}
\begin{claim}
	The construction of $A_\alpha$, $\alpha<\continuum$ as above is possible.
\end{claim}
\prfofClaim
Suppose that $\alpha\in\continuum\setminus\omega$ and $A_\beta$,
$\beta<\alpha$ have been constructed according to \xitemof{f-0}, \xitemof{f-1} 
and \xitemof{f-2}. Let 
\begin{xitemize}
\item[]
	$S_\alpha=\setof{\beta<\alpha}{
           B_{\beta}\notin \ide(\setof{A_\delta}{{\delta}<{\alpha}})}
          $. 
\end{xitemize}
Let
	$\poP_\alpha=\setof{\pairof{\varphi,s}}{\varphi\in\Fn(T,2),\,s\in[\alpha]^{<\aleph_0}}$
be the \po\ with the ordering defined by
\begin{xitemize}
\item[] 
	$\pairof{\varphi',s'}\leq_{\poP_\alpha}\pairof{\varphi,s}$\ \ $\Leftrightarrow$\\[\jot]
	\phantom{$\pairof{\varphi',s'}\leq$}%
	$\varphi\subseteq \varphi'$, $s\subseteq s'$ and\\
	\phantom{$\pairof{\varphi',s'}\leq$}%
	$\forall t\in\dom(\varphi')\setminus\dom(\varphi)\ 
	(\varphi'(t)=1\ \rightarrow\ t\not\in A_\delta\mbox{ for all }\delta\in s)$
\end{xitemize}
for $\pairof{\varphi,s}$, $\pairof{\varphi',s'}\in\poP_\alpha$. 

$\poP_\alpha$ is $\sigma$-centered since $\pairof{\varphi,s}$,
$\pairof{\varphi',s'}\in\poP_\alpha$ are compatible if $\varphi=\varphi'$.  

For $\beta<\alpha$, let 
\begin{xitemize}
\item[] $C_\beta=\setof{\pairof{\varphi,s}\in\poP_\alpha}{\beta\in s}$
\end{xitemize}
and, for $\beta\in S_\alpha$ and $n\in\omega$, let 
\begin{xitemize}
\item[] $D_{\beta,n}=\setof{\pairof{\varphi,s}\in\poP_\alpha}{
	\exists t\in\dom(\varphi)\ (\ell(t)\geq n\ \land\ \varphi(t)=1\ \land\ t\in 
	B_{\beta})}$. 
\end{xitemize}
It is easy to see that $C_\beta$, $\beta<\alpha$ and $D_{\beta,n}$,
$\beta\in S_\alpha$, $n\in\omega$ are dense in $\poP_\alpha$. 
Let
\begin{xitemize}
\item[] $\calD=\setof{C_{\beta}}{\beta<\alpha}
	\cup\setof{D_{\beta,n}}{\beta\in S_\alpha,\,n\in\omega}
	$.
\end{xitemize}
Since $\cardof{\calD}<\continuum$, we can 
apply $\MA(\sigma\mbox{-centered})$ to obtain 
a $(\calD,\poP_\alpha)$-generic filter $G$. Let 
\begin{xitemize}
\item[] $A_\alpha=\setof{t\in T}{\varphi(t)=1\mbox{ for some }\pairof{\varphi,s}\in G}$.
\end{xitemize}
Then this $A_\alpha$ is as desired. 
\qedofClaim\qedskip

Let $\calF=\setof{A_\alpha}{\alpha<\continuum}$. $\calF$ is infinite by 
\xitemof{f-1} and mad  by \xitemof{f-2}.

We show that $\calF$ is $\continuum$-almost decided. First, note that we have 
$\gma=\continuum$ by the assumptions of the theorem. 
By \xitemof{f-2}, we have: 
\begin{xitemize}
	\xitem[f-4] For any $B\in[T]^{\aleph_0}$, if  
$B\notin \ide(\setof{A_{\alpha}}{{\alpha}<\continuum})$ then \\
$\cardof{\setof{{\alpha}<\continuum}{\cardof{A_{\alpha}\cap B}<\aleph_0}}<\continuum$.
\end{xitemize}
Suppose that $\calH\in[\calF]^\continuum$ and 
$B\in\calH^\perp$. Then $\cardof{\setof{{\alpha}<\continuum}{\cardof{A_{\alpha}\cap B}<\aleph_0}}=\continuum$ and so $B\in \ide(\calF)$ by \xitemof{f-4}.
Thus there is a finite $\calF'\subset \calF$ such that 
$B\subset^* \cup\calF'$ and $F\cap B$ is infinite for each $F\in \calF'$.
But $B\in \calH^\perp$ so $\calF'\cap \calH=\emptyset$. Thus
$\calF'$ witnesses that $B\in \ide(\calF\setminus \calH)$
which was to be proved.
\qed
\end{proof}

For a mad family $\calF$ on $T$, $\calC\subseteq\calF$ is said to be 
{\em minimal in $\calF$} if $\gma^+(\calF\setminus\calC)=\cardof{\calC}$. A mad 
family $\calF$ is said to be {\em$\lambda$-minimal\/} if every 
$\calC\in[\calF]^{\lambda}$ is minimal in $\calF$. 

\begin{lemma}\label{almost=decided-minimal}
	Suppose that $\calF$ is a mad family on $T$. 

\assert{1} If $\calF$ is $\cardof{\calF}$-minimal then $\cardof{\calF}=\gma$.

\assert{2} If $\calB\subseteq\calF$ almost decides $\calF$ and
$\calF\setminus\calB$ is infinite then  
$\calF\setminus\calB$ is minimal in $\calF$. 

\assert{3} If $\calF$ is $\kappa$-almost decided for $\kappa=\cardof{\calF}$ 
then $\calF$ is $\lambda$-minimal for all $\omega\leq\lambda<\kappa$. 

\assert{4} If $\cardof{\calF}=\gma$ and $\calF$ is $\gma$-almost decided then 
$\calF$ is $\gma$-minimal. 
\end{lemma}
\begin{proof} \assertof{1}: If $\calF$ is $\cardof{\calF}$-minimal then 
$\calF$ itself is minimal in $\calF$. Thus   
$\gma=\gma^+(\emptyset)=\gma^+(\calF\setminus\calF)=\cardof{\calF}$.  
\smallskip

\assertof{2}: First, note that, for any infinite ad $\calF$, we have
$\gma(\ide(\calF))=\cardof{\calF}$. 

Suppose that $\calF$ is a mad family on $T$ and $\calB\subseteq\calF$ 
almost decides $\calF$, i.e.\ $\calB^\perp=\ide(\calF\setminus\calB)$. 
Hence 
\begin{xitemize}
\item[] $\gma^+(\calF\setminus(\calF\setminus\calB))=\gma^+(\calB)
	=\gma(\calB^\perp)=\gma(\ide(\calF\setminus\calB))=\cardof{\calF\setminus\calB}$. 
\end{xitemize}

\assertof{3}: Suppose that $\kappa=\cardof{\calF}$ and $\calF$ 
is $\kappa$-almost decided. If $\calC\in[\calF]^{\lambda}$ for some 
$\omega\leq\lambda<\kappa$ then $\cardof{\calF\setminus\calC}=\kappa$ and 
hence $\calF\setminus\calC$ almost decides $\calF$. By \assertof{2} it follows 
that $\calC=\calF\setminus(\calF\setminus\calC)$ is minimal in $\calF$. 

\assertof{4}: Suppose that $\cardof{\calF}=\gma$ and $\calF$ 
is $\gma$-almost decided. Suppose that $\calC\in[\calF]^\gma$. If 
$\cardof{\calF\setminus\calC}<\gma$, then clearly
$\gma^+(\calF\setminus\calC)=\gma=\cardof{\calC}$. Hence $\calC$ is minimal 
in $\calF$.  
If $\cardof{\calF\setminus\calC}=\gma$ then $\calF\setminus\calC$  
almost decides $\calF$. Thus, by \assertof{2}, 
$\calC=\calF\setminus(\calF\setminus\calC)$ is again minimal in $\calF$. 
\qed
\end{proof}
\begin{corollary}
	\label{c-minimal}
	Assume that $\MA(\sigma\mbox{-centered\/{}})$ holds. 
	Then there is a mad 
	family $\calF$ on $T$ which is $\lambda$-minimal for all
	$\omega\leq\lambda\leq\continuum$.  
\end{corollary}
\begin{proof}
By \Thmof{c-almost decided} and 
\Lemmaof{almost=decided-minimal},\,\assertof{3}, \assertof{4}. 
\qed
\end{proof}
\Thmof{c-almost decided} can be further improved to the following theorem:
\begin{theorem}
	\label{c-almost=decided-x}
	Assume that $\MA(\sigma\mbox{-centered\/{}})$ holds. Let $\kappa=\continuum$.
	Then there is a $\Cohen{\omega}$-indestructible mad family $\calF$ 
	(of size $\kappa$) \st\ 
	\begin{xitemize}
	  \xitem[k-a-d-0] $V^\Cohen{\omega}\models\calF\mbox{ is }\kappa\mbox{-almost 
		decided on }T$.
	\end{xitemize}
\end{theorem}
\begin{proof} Let $\seqof{\pairof{t_\beta,\utilde{B}_\beta}}{\beta<\kappa}$ be an 
enumeration of 
\begin{xitemize}
\item[] $T\times\setof{\utilde{B}}{\utilde{B}
	\mbox{ is a nice }\Cohen{\omega}\mbox{-name of an element of }
			 [T]^{\aleph_0}\mbox{ in }V^{\Cohen{\omega}}}$. 
\end{xitemize}
Let $A_\alpha$, $\alpha<\kappa$ be then defined inductively just as in the 
proof of \Thmof{c-almost decided} with 
\begin{xitemize}
	\item[\xitemof{f-2}$'$] For $\beta<\alpha$, if
	$t\forces{\Cohen{\omega}}{
		\utilde{B}_{\alpha}\notin \ide(\setof{A_\delta}{{\delta}<{\alpha}})}$ 
	then
	$t\forces{\Cohen{\omega}}{\cardof{A_\alpha\cap \utilde{B}_\beta}=\aleph_0}$ 
\end{xitemize}
in place of \xitemof{f-2}. 
\qed
\end{proof}
\begin{corollary}
	\label{c-almost=decided-minimal-x}
  For any cardinal $\kappa\geq\continuum$ in the ground model $V$ there is a 
  cardinal preserving generic extension $W$ of $V$ \st, in $W$, $\kappa<\continuum$
  and 
	there is a $\kappa$-almost decided mad family $\calF$ of size $\kappa$ 
	(furthermore $\calF$ is $\lambda$-minimal for all
	$\omega\leq\lambda\leq\kappa$).  
\end{corollary}
\begin{proof}
  First extend $V$ to a model $V'$ of $\kappa=\continuum$ and $\MA(\sigma\mbox{-centered})$. 
In $V'$, let $\calF$ be as in \Thmof{c-almost=decided-x}. Then $\calF$ is as 
desired in $V^\Cohen{\mu}$ for any $\mu>\kappa$. The claim in the 
parentheses follows from 
\Lemmaof{almost=decided-minimal},\,\assertof{3} and \xitemof{f-2}$'$.  
\qed
\end{proof}

%
%

\end{document}